\documentclass{article}
\usepackage[utf8]{inputenc}
\usepackage{cite}
\usepackage[english]{babel}
\usepackage{subfig}
\usepackage{amsmath}
\usepackage{amssymb}
\usepackage{amsthm}
\usepackage{tasks}
\usepackage[table]{xcolor}
\usepackage[inline]{enumitem}
\usepackage{mathrsfs}
\usepackage{bm}
\usepackage{dsfont}
\usepackage{fontawesome}

\makeatletter
\newcommand{\thickhline}{%
    \noalign {\ifnum 0=`}\fi \hrule height 1pt
    \futurelet \reserved@a \@xhline
}
\newcolumntype{"}{@{\hskip\tabcolsep\vrule width 1pt\hskip\tabcolsep}}
\makeatother

\usepackage{empheq}
\usepackage{float}
\usepackage{booktabs,hyperref}
\usepackage{float,caption,hypcap}
\newtheorem{theorem}{Theorem}[section]
\newtheorem{lemma}[theorem]{Lemma}

\title{Regular Trace Formula of Eigenvalues of a Discontinuous One-point Boundary Value Problem  with Retarded Argument}

\usepackage{tikz,xcolor,hyperref}
\usepackage{authblk}
\definecolor{lime}{HTML}{A6CE39}
\DeclareRobustCommand{\orcidicon}{
	\begin{tikzpicture}
	\draw[lime, fill=lime] (0,0)
	circle [radius=0.16]
	node[white] {{\fontfamily{qag}\selectfont \tiny ID}};
	\draw[white, fill=white] (-0.0625,0.095)
	circle [radius=0.007];
	\end{tikzpicture}
	\hspace{-2mm}}
\foreach \x in {A, ..., Z}{\expandafter\xdef\csname orcid\x\endcsname{\noexpand\href{https://orcid.org/\csname orcidauthor\x\endcsname}
			{\noexpand\orcidicon}}
}

\author[1]{Yunus SAÇLI$^1$\orcidA{}, Seda KIZILBUDAK ÇALIŞKAN$^2$\orcidB{}\\
$^1$ Ph.D Student at Institute of Natural Science, Yildiz Technical University, 34210, Istanbul, Turkey\\
yunus.sacli@std.yildiz.edu.tr\\
$^2$Department of Mathematics, Faculty of Arts and Science, Yildiz Technical University, 34210, Istanbul, Turkey\\
skizilb@yildiz.edu.tr}

\begin{document}

\maketitle

\begin{abstract} 
In this study, we found a regular trace formula for the eigenvalues of the boundary value problem, which we created with the second-order differential equation with eigen parameter and discontinuity at  $x={\pi \over 2}$, which is an interior point of the finite range $[0,\pi]$, and boundary conditions that also contain eigen parameter, and interface conditions.
\end{abstract}

\maketitle

\section{Introduction}
\noindent The purpose of this study is to obtain the regular trace formula for the  eigenvalues of the boundary value problem formed by the second order differential equation contains a eigen parameter
\begin{equation}
    y''(x) + q(x)y(x - \Delta (x)) =  - {\mu ^2}y(x)
    \end{equation}

\noindent on   $[0,{\textstyle{\pi  \over 2}}) \cup ({\textstyle{\pi  \over 2}},\pi ]$ , with  boundary conditions depend on a eigen parameter

\begin{equation}
(\mu {a'_1} + {a_1})y(0) - (\mu {a'_2} + {a_2})y'(0) = 0 
\end{equation}

\begin{equation}
    \cos by(\pi ) + \mu \sin by'(\pi ) = 0
\end{equation}  
\noindent and with interface conditions
\begin{equation}
    y({\textstyle{\pi  \over 2}} - 0) = \delta y({\textstyle{\pi  \over 2}} + 0)
\end{equation}
\begin{equation}
   y'({\textstyle{\pi  \over 2}} - 0) = \delta y'({\textstyle{\pi  \over 2}} + 0).
\end{equation}
\noindent  $q(x)$ is the real-valued continuous function  on both the intervals $[0,{\textstyle{\pi  \over 2}})$ and $({\textstyle{\pi  \over 2}},\pi ]$   and $q(x)$ has finite limits

\begin{equation}
\left\{\begin{array}{rl}
 q({\textstyle{\pi  \over 2}} + 0) =    &   \mathop {\lim }\limits_{x \to {\textstyle{\pi  \over 2}} + 0} q(x)\\
  q({\textstyle{\pi  \over 2}} - 0) =   &  \mathop {\lim }\limits_{x \to {\textstyle{\pi  \over 2}} - 0} q(x)
\end{array} \right.
\end{equation}

\noindent and $\Delta (x) \ge 0$ is the real-values function and continuous on $[0,{\textstyle{\pi  \over 2}}) \cup ({\textstyle{\pi  \over 2}},\pi ]$ and has finite limits

\begin{equation}
\left\{\begin{array}{rl}
\Delta ({\textstyle{\pi  \over 2}} + 0) =     &  \mathop {\lim }\limits_{x \to {\textstyle{\pi  \over 2}} + 0} \Delta (x) \\
\Delta ({\textstyle{\pi  \over 2}} - 0) =     & \mathop {\lim }\limits_{x \to {\textstyle{\pi  \over 2}} - 0} \Delta (x),
\end{array} \right.
\end{equation}

 Moreover,
\begin{equation*}
\left\{\begin{array}{rl}
 x - \Delta (x) \ge 0 &\mbox{if}\quad x \in \left[ {0,{\textstyle{\pi  \over 2}}} \right)\\
 x - \Delta (x) \ge {\textstyle{\pi  \over 2}}&\mbox{if}\quad x \in \left( {{\textstyle{\pi  \over 2}},\pi } \right]\end{array}\right.
\end{equation*}

\noindent $\mu $ is a eigen parameter,  ${a'_1},{a_1},{a'_2},{a_2},\delta  \ne 0$ are arbitrary real numbers.

\noindent The scholars have long been interested in the theory of regular trace of ordinary differential operators. Gelfand and Levitan \cite{gelfand} firstly obtained the trace formula for the Sturm-Liouville differential equation. After this study mathematicians were interested in developing trace formulas for different differential operators. In the survey article \cite{sadovnichii}, the whole history of the regular trace of the linear operators is given. The regular trace formulas for differential equations are found \cite{boglu86, sahinturk, makin, maksudov, orucoglu} and however, there are a small number of works on the regular trace for the differential equations with retarded argument.

\noindent The boundary value problems with a discontinuity condition inside the interval have many applications such as
mathematics, mechanics, radio electronics, geophysics and other fields of science. Furthermore, \cite{norkin} contains several physical applications of such problems. Boundary value problems involving interface conditions are commonly encountered in the theory of heat and mass transfer in a wide range of physical transfer problems.

\noindent Due to the demands of modern technology, engineering, and physics, the problems with interface conditions have grown in importance in recent years. The retarded differential equations, which are one of the differential classes, can be continuous in the interval in they are defined or have discontinuity at one or more interior point of this interval. This type of problem has been addressed in many studies \cite{yang, norkin, hira, sen, bayramoglu, caliskan}.

\noindent In our study, we created our problem by adding the eigen parameter-containing conditions to our equation with a discontinuity at one point, as well as the interface conditions in order not to complicate the solution of discontinuity-related irregularities, and we have obtained a regular trace formula for the eigenvalues of this problem
\[\begin{array}{l}
\mu _{ - 0}^2 + \mu _0^2 + \sum\limits_{n = 1}^\infty  {\left( {\mu _{ - n}^2 + \mu _n^2 - 2{{(n - 1)}^2} + \frac{4}{\pi }\left( { - \frac{{{{a'_1}}}}{{{{a'_2}}}} + P(\mu ,\Delta (\theta )) + Q(\mu ,\Delta (\theta ))} \right)} \right)} \\
\,\,\,\,\,\,\, =  - \frac{2}{\pi }\left( { - \frac{{{{a'_1}}}}{{{{a'_2}}}} + P(0,\Delta (\theta )) + Q(0,\Delta (\theta ))} \right)\\
\,\,\,\,\,\,\,\,\,\,\, - {\left( { - \frac{{{{a'_1}}}}{{{{a'_2}}}} + P(0,\Delta (\theta )) + Q(0,\Delta (\theta ))} \right)^2}\\
\,\,\,\,\,\,\,\,\,\,\, + {\left( {\frac{{{a_2}}}{{{{a'_2}}}} + R(0,\Delta (\theta )) + S(0,\Delta (\theta ))} \right)^2} + O\left( {\frac{1}{{{N_0}}}} \right).
\end{array}\]

\noindent We think that this study will make new contributions to the problems in the field of quantum statistics and technology.

\section{Transform of Boundary Value Problem to Integral Equation}

\noindent Let ${\omega _1}(x,\mu )$ and ${\omega _2}(x,\mu )$ be the solutions of our quadratic equation with discontinuity at $x = {\textstyle{\pi  \over 2}}$. ${\omega _1}(x,\mu )$ and  ${\omega ' _1}(x,\mu )$ satisfy the following conditions 
\begin{equation}
\left\{\begin{array}{rl}
      {\omega _1}(0,\mu ) =   & \mu {a'_2} + {a_2} \\
  {\omega '_1}(0,\mu ) =   & \mu {a'_1} + {a_1}
\end{array} \right.
\end{equation}

\noindent on $[0,{\textstyle{\pi  \over 2}}]$. Under these conditions, ${\omega _1}(x,\mu )$ is the unique solution of our equation on $[0,{\textstyle{\pi  \over 2}}]$.

\noindent We can express ${\omega _2}(x,\mu )$ which is the solution of our equation, using interface conditions, depending on ${\omega _1}(x,\mu )$ in the interval $[{\textstyle{\pi  \over 2}},\pi ]$ as follows:
\begin{equation}
\left\{\begin{array}{rl}
     {\omega _2}( \frac{\pi}{2},\mu ) =  &   {\delta ^{ - 1}}{\omega _1}( \frac{\pi}{2},\mu ) \\
  {\omega '_2}( \frac{\pi}{2},\mu ) =   &  {\delta ^{ - 1}}{\omega '_1}( \frac{\pi}{2},\mu ).
\end{array} \right.
\end{equation}

\noindent With these conditions, ${\omega _2}(x,\mu )$ is the unique solution of our equation on the interval $[{\textstyle{\pi  \over 2}},\pi ]$.

\noindent As a result, the unique solution of our equation defined in the interval $[0,{\textstyle{\pi  \over 2}}] \cup [{\textstyle{\pi  \over 2}},\pi ]$ which satisfies the first of the boundary condition and interface conditions, is as follows
\begin{equation}
    \omega (x,\mu ) = \left\{ \begin{array}{l}
{\omega _1}(x,\mu ),\,\,\,x \in [0,{\textstyle{\pi  \over 2}})\\
{\omega _2}(x,\mu ),\,\,\,x \in ({\textstyle{\pi  \over 2}},\pi ].
\end{array} \right.
\end{equation}

\noindent The solution ${\omega }(x,\mu )$  defined above is nontrivial solution of the second order differential equation $(1)$ that satisfies boundary condition $(2)$ and interface conditions $(4)-(5)$.

\begin{lemma} Let $\omega (x,\mu )$ be a solution of the second order differential equation $(1)$ and $\mu  > 0$. Thus providing the following integral equations :

\begin{equation}
\begin{array}{rl}
 {\omega _1}(x,\mu ) =     & 
    \cos \mu x(\mu {a'_2} + {a_2}) + \frac{{\sin \mu x}}{\mu }(\mu {a'_1} + {a_1}) \\
     &  - \frac{1}{\mu }\int\limits_0^x {q(\theta )\sin \mu } (x - \theta ){\omega _1}(\theta  - \Delta (\theta ),\mu )d\theta
\end{array}
\end{equation}
and
\begin{equation}
  \begin{array}{l}
{\omega _2}(x,\mu ) = \frac{1}{{\mu \delta }}\sin \mu (x - {\raise0.5ex\hbox{$\scriptstyle \pi $}
\kern-0.1em/\kern-0.15em
\lower0.25ex\hbox{$\scriptstyle 2$}}){\omega '_1}({\raise0.5ex\hbox{$\scriptstyle \pi $}
\kern-0.1em/\kern-0.15em
\lower0.25ex\hbox{$\scriptstyle 2$}},\mu )\\
\,\,\,\,\,\,\,\,\,\,\,\,\,\,\,\,\,\,\,\,\,\,\,\,\,\,\,+ \frac{1}{\delta }\cos \mu (x - {\raise0.5ex\hbox{$\scriptstyle \pi $}
\kern-0.1em/\kern-0.15em
\lower0.25ex\hbox{$\scriptstyle 2$}}){\omega _1}({\raise0.5ex\hbox{$\scriptstyle \pi $}
\kern-0.1em/\kern-0.15em
\lower0.25ex\hbox{$\scriptstyle 2$}},\mu )\\
\,\,\,\,\,\,\,\,\,\,\,\,\,\,\,\,\,\,\,\,\,\,\,\,\,\,\, - \frac{1}{\mu }\int\limits_{{\raise0.5ex\hbox{$\scriptstyle \pi $}
\kern-0.1em/\kern-0.15em
\lower0.25ex\hbox{$\scriptstyle 2$}}}^x {q(\theta )} \sin \mu (x - \theta ){\omega _2}(\theta  - \Delta (\theta ),\mu )d\theta. 
\end{array}
\end{equation}
\end{lemma}
\begin{proof} To prove this, it is enough to substitute $ - {\omega ''_1}(\theta ,\mu ) - {\mu ^2}{\omega _1}(\theta ,\mu )$ and \linebreak $ - {\omega ''_2}(\theta ,\mu ) - {\mu ^2}{\omega _2}(\theta ,\mu )$ instead of $q(\theta ){\omega _1}(\theta  - \Delta (\theta ),\mu )$ and $q(\theta ){\omega _2}(\theta  - \Delta (\theta ),\mu )$ in integrals in equations $(11)$ and $(12)$, respectively, and integrate by parts twice.
\end{proof}

\noindent Let's call the equation we obtained by considering $\omega (x,\mu )$ in $(2)$ as the characteristic equation
\begin{equation}
    F(\mu ) = \cos b\omega (\pi ,\mu ) + \mu \sin b\omega '(\pi ,\mu ) = 0.
\end{equation}

\noindent The eigenvalues of our problem $(1)-(5)$ coincides with the roots of $F(\mu )$, and the eigenvalues are simple.

\noindent $F(\mu )$ is determined with equations $(11)$, $(12)$ and their derivatives.
\begin{equation}
\begin{array}{rl}
  {\omega '_1}(x,\mu ) =    & - \mu \sin \mu x(\mu {a'_2} + {a_2}) + \cos \mu x(\mu {a'_1} + {a_1})\\
     &  - \int\limits_0^x {q(\theta )\cos \mu } (x - \theta ){\omega _1}(\theta  - \Delta (\theta ),\mu )d\theta 
\end{array}
\end{equation}

\noindent and

\begin{equation}
   \begin{array}{l}
{{\omega '_2}}(x,\mu ) = \frac{1}{\delta }\cos \mu (x - {\raise0.5ex\hbox{$\scriptstyle \pi $}
\kern-0.1em/\kern-0.15em
\lower0.25ex\hbox{$\scriptstyle 2$}}){{\omega '}_1}({\raise0.5ex\hbox{$\scriptstyle \pi $}
\kern-0.1em/\kern-0.15em
\lower0.25ex\hbox{$\scriptstyle 2$}},\mu ) \\
\,\,\,\,\,\,\,\,\,\,\,\,\,\,\,\,\,\,\,\,\,\,\,\,\,\,\,\,- \frac{\mu }{\delta }\sin \mu (x - {\raise0.5ex\hbox{$\scriptstyle \pi $}
\kern-0.1em/\kern-0.15em
\lower0.25ex\hbox{$\scriptstyle 2$}}){\omega _1}({\raise0.5ex\hbox{$\scriptstyle \pi $}
\kern-0.1em/\kern-0.15em
\lower0.25ex\hbox{$\scriptstyle 2$}},\mu )\\
\,\,\,\,\,\,\,\,\,\,\,\,\,\,\,\,\,\,\,\,\,\,\,\,\,\,\,\, - \int\limits_{{\raise0.5ex\hbox{$\scriptstyle \pi $}
\kern-0.1em/\kern-0.15em
\lower0.25ex\hbox{$\scriptstyle 2$}}}^x {q(\theta )} \cos \mu (x - \theta ){\omega _2}(\theta  - \Delta (\theta ),\mu )d\theta.
\end{array}
\end{equation}

\noindent Then we have:

\begin{equation*}
\begin{array}{rl}
F(\mu ) =& \frac{{\left( {\mu {{a'_1}} + {a_1}} \right)}}{{\mu \delta }}\cos b .\sin \mu \pi  + \frac{{\left( {\mu {{a'_2}} + {a_2}} \right)}}{\delta }\cos b .\cos \mu \pi \\
& - \frac{{\cos b }}{{\mu \delta }}\int\limits_0^{{\raise0.5ex\hbox{$\scriptstyle \pi $}
\kern-0.1em/\kern-0.15em
\lower0.25ex\hbox{$\scriptstyle 2$}}} {q(\theta )} \sin \mu (\pi  - \theta ){\omega _1}(\theta  - \Delta (\theta ),\mu )d\theta \\
& - \frac{{\cos b }}{\mu }\int\limits_{{\raise0.5ex\hbox{$\scriptstyle \pi $}
\kern-0.1em/\kern-0.15em
\lower0.25ex\hbox{$\scriptstyle 2$}}}^\pi  {q(\theta )} \sin \mu (\pi  - \theta ){\omega _2}(\theta  - \Delta (\theta ),\mu )d\theta \\
& + \frac{{\mu \left( {\mu {{a'_1}} + {a_1}} \right)}}{\delta }\sin b .\cos \mu \pi  + \frac{{{\mu ^2}\left( {\mu {{a'_2}} + {a_2}} \right)}}{\delta }\sin b .\sin \mu \pi \\
& - \frac{{\mu \sin b }}{\delta }\int\limits_0^{{\raise0.5ex\hbox{$\scriptstyle \pi $}
\kern-0.1em/\kern-0.15em
\lower0.25ex\hbox{$\scriptstyle 2$}}} {q(\theta )} \cos \mu (\pi  - \theta ){\omega _1}(\theta  - \Delta (\theta ),\mu )d\theta \\
& - \mu \sin b \int\limits_{{\raise0.5ex\hbox{$\scriptstyle \pi $}
\kern-0.1em/\kern-0.15em
\lower0.25ex\hbox{$\scriptstyle 2$}}}^\pi  {q(\theta )} \cos \mu (\pi  - \theta ){\omega _2}(\theta  - \Delta (\theta ),\mu )d\theta  \\
=& 0.
\end{array}    
\end{equation*}

\noindent Let's use the sequential approximation of ${\omega _1}(x,\mu )$ and ${\omega _2}(x,\mu )$ in this last equation of $F(\mu )$. For this we write the sequential approximation:

\begin{equation}
    \begin{array}{rl}
{\omega _1}(x, \mu ) &= \frac{{\sin  \mu x}}{ \mu }( \mu {{a'_1}} + {a_1}) + \cos  \mu x( \mu {{a'_2}} + {a_2}) \\
&+ \frac{{ \mu {{a'_1}} + {a_1}}}{{2{ \mu ^2}}}\int\limits_0^x {q(\theta )\cos  \mu (x - \Delta (\theta } ))d\theta \\
& - \frac{{ \mu {{a'_1}} + {a_1}}}{{2{ \mu ^2}}}\int\limits_0^x {q(\theta )\cos  \mu (x - 2\theta  + \Delta (\theta ))d\theta } \\
& - \frac{{ \mu {{a'_2}} + {a_2}}}{{2 \mu }}\int\limits_0^x {q(\theta )\sin  \mu (x - \Delta (\theta ))} d\theta \\
& - \frac{{ \mu {{a'_2}} + {a_2}}}{{2 \mu }}\int\limits_0^x {q(\theta )\sin  \mu (x - 2\theta  + \Delta (\theta ))d\theta. } 
\end{array}
\end{equation}

\begin{equation}
    \begin{array}{rl}
{\omega _2}(x, \mu ) =&  \frac{{\sin  \mu x}}{{ \mu \delta }}( \mu {{a'}_1} + {a_1}) + \frac{{\cos  \mu x}}{\delta }( \mu {{a'}_2} + {a_2}) \\
& + \frac{{ \mu {{a'}_1} + {a_1}}}{{2{ \mu ^2}\delta }}\int\limits_0^x {q(\theta )\cos  \mu (x - \Delta (\theta } ))d\theta \\
&  - \frac{{ \mu {{a'}_1} + {a_1}}}{{2{ \mu ^2}\delta }}\int\limits_0^x {q(\theta )\cos  \mu (x - 2\theta  + \Delta (\theta ))d\theta } \\
&  - \frac{{ \mu {{a'}_2} + {a_2}}}{{2 \mu \delta }}\int\limits_0^x {q(\theta )\sin  \mu (x - \Delta (\theta ))} d\theta \\
&  - \frac{{ \mu {{a'}_2} + {a_2}}}{{2 \mu \delta }}\int\limits_0^x {q(\theta )\sin  \mu (x - 2\theta  + \Delta (\theta ))d\theta. } 
\end{array}
\end{equation}

\noindent Then we have,

\[\begin{array}{rl}
F( \mu ) &=  \frac{{( \mu {{a'_1}} + {a_1})}}{{ \mu \delta }}\cos b .\sin  \mu \pi  + \frac{{( \mu {{a'_2}} + {a_2})}}{\delta }\cos b .\cos  \mu \pi \\
&  - \frac{{\cos b ( \mu {{a'_1}} + {a_1})}}{{{ \mu ^2}\delta }}\int\limits_0^\pi  {q(\theta )} \sin  \mu (\pi  - \theta )\sin  \mu (\theta  - \Delta (\theta ))d\theta \\
&  - \frac{{\cos b ( \mu {{a'_2}} + {a_2})}}{{ \mu \delta }}\int\limits_0^\pi  {q(\theta )} \sin  \mu (\pi  - \theta )\cos  \mu (\theta  - \Delta (\theta ))d\theta \\
&  - \frac{{\cos b ( \mu {{a'_1}} + {a_1})}}{{2{ \mu ^3}\delta }}\int\limits_0^\pi  {q(\theta )} \sin  \mu (\pi  - \theta )\int\limits_0^{\theta  - \Delta (\theta )} {q({t_1})} \cos  \mu (\theta  - \Delta (\theta ) - \Delta ({t_1}))d{t_1}d\theta \\
&  - \frac{{\cos b ( \mu {{a'_1}} + {a_1})}}{{2{ \mu ^3}\delta }}\int\limits_0^\pi  {q(\theta )} \sin  \mu (\pi  - \theta )\int\limits_0^{\theta  - \Delta (\theta )} {q({t_1})} \cos  \mu (\theta  - \Delta (\theta ) - 2{t_1} + \Delta ({t_1}))d{t_1}d\theta \\
&  + \frac{{\cos b ( \mu {{a'_2}} + {a_2})}}{{2{ \mu ^2}\delta }}\int\limits_0^\pi  {q(\theta )} \sin  \mu (\pi  - \theta )\int\limits_0^{\theta  - \Delta (\theta )} {q({t_1})} \sin  \mu (\theta  - \Delta (\theta ) - \Delta ({t_1}))d{t_1}d\theta \\
&  + \frac{{\cos b ( \mu {{a'_2}} + {a_2})}}{{2{ \mu ^2}\delta }}\int\limits_0^\pi  {q(\theta )} \sin  \mu (\pi  - \theta )\int\limits_0^{\theta  - \Delta (\theta )} {q({t_1})} \sin  \mu (\theta  - \Delta (\theta ) - 2{t_1} + \Delta ({t_1}))d{t_1}d\theta \\
&  + \frac{{ \mu ( \mu {{a'_1}} + {a_1})}}{\delta }\sin b .\cos  \mu \pi  - \frac{{{ \mu ^2}( \mu {{a'_2}} + {a_2})}}{\delta }\sin b .\sin  \mu \pi \\
&  - \frac{{\sin b ( \mu {{a'_1}} + {a_1})}}{\delta }\int\limits_0^\pi  {q(\theta )} \cos  \mu (\pi  - \theta )\sin  \mu (\theta  - \Delta (\theta ))d\theta \\
&  - \frac{{ \mu ( \mu {{a'_2}} + {a_2})\sin b }}{\delta }\int\limits_0^\pi  {q(\theta )} \cos  \mu (\pi  - \theta )\cos  \mu (\theta  - \Delta (\theta ))d\theta \\
&  - \frac{{\sin b ( \mu {{a'_1}} + {a_1})}}{{2 \mu \delta }}\int\limits_0^\pi  {q(\theta )} \cos  \mu (\pi  - \theta )\int\limits_0^{\theta  - \Delta (\theta )} {q({t_1})} \cos  \mu (\theta  - \Delta (\theta ) - \Delta ({t_1}))d{t_1}d\theta \\
&  + \frac{{\sin b ( \mu {{a'_1}} + {a_1})}}{{2 \mu \delta }}\int\limits_0^\pi  {q(\theta )} \cos  \mu (\pi  - \theta )\int\limits_0^{\theta  - \Delta (\theta )} {q({t_1})} \cos  \mu (\theta  - \Delta (\theta ) - 2{t_1} + \Delta ({t_1}))d{t_1}d\theta \end{array}\]
\[\begin{array}{rl}&  \,\,\,\,\,\,\,\,\,\,\,\,+ \frac{{\sin b ( \mu {{a'_2}} + {a_2})}}{{2\delta }}\int\limits_0^\pi  {q(\theta )} \cos  \mu (\pi  - \theta )\int\limits_0^{\theta  - \Delta (\theta )} {q({t_1})} \sin  \mu (\theta  - \Delta (\theta ) - \Delta ({t_1}))d{t_1}d\theta \\
&  \,\,\,\,\,\,\,\,\,\,\,\,+ \frac{{\sin b ( \mu {{a'_2}} + {a_2})}}{{2\delta }}\int\limits_0^\pi  {q(\theta )} \cos  \mu (\pi  - \theta )\int\limits_0^{\theta  - \Delta (\theta )} {q({t_1})} \sin  \mu (\theta  - \Delta (\theta ) - 2{t_1} + \Delta ({t_1}))d{t_1}d\theta  \\
&\,\,\,\,\,\,\,\,\,\,\,\,= 0.
\end{array}\]

\noindent By making trigonometric transformations and some adjustments, we find:

\[\begin{array}{rl}
F ( \mu ) &=   - \frac{{{ \mu ^3}{{a'_2}}\sin b}}{\delta }\sin  \mu \pi  + \frac{{ \mu {{a'_2}}}}{\delta }\cos b\cos  \mu \pi  + \frac{{ \mu {a_1}}}{\delta }\sin b\cos  \mu \pi \\
&  + \frac{{ \mu {a_1}\sin b}}{{2\delta }}\int\limits_0^\pi  {q(\theta )} \sin  \mu \Delta (\theta )\cos  \mu \pi d\theta  \\
&- \frac{{ \mu {{a'_1}}\sin b}}{{2\delta }}\int\limits_0^\pi  {q(\theta )} \sin  \mu (2\theta  - \Delta (\theta ))\cos  \mu \pi d\theta \\
&  - \frac{{ \mu {a_2}\sin b}}{{2\delta }}\int\limits_0^\pi  {q(\theta )} \cos  \mu \Delta (\theta )\cos  \mu \pi d\theta  \\
&- \frac{{ \mu {a_2}\sin b}}{{2\delta }}\int\limits_0^\pi  {q(\theta )} \cos  \mu (2\theta  - \Delta (\theta ))\cos  \mu \pi d\theta \\
&  - \frac{{ \mu {{a'_1}}\sin b}}{{2\delta }}\int\limits_0^\pi  {q(\theta )} \cos  \mu \Delta (\theta )\sin  \mu \pi d\theta \\
& - \frac{{ \mu {{a'_1}}\sin b}}{{2\delta }}\int\limits_0^\pi  {q(\theta )} \cos  \mu (2\theta  - \Delta (\theta ))\sin  \mu \pi d\theta \\
&  - \frac{{ \mu {a_2}\sin b}}{{2\delta }}\int\limits_0^\pi  {q(\theta )} \sin  \mu \Delta (\theta )\sin  \mu \pi d\theta  \\
&- \frac{{ \mu {a_2}\sin b}}{{2\delta }}\int\limits_0^\pi  {q(\theta )} \sin  \mu (2\theta  - \Delta (\theta ))\sin  \mu \pi d\theta \\
&  + \frac{{{ \mu ^2}{{a'_1}}\sin b}}{\delta }\sin  \mu \pi  - \frac{{{ \mu ^2}{{a'_1}}\sin b}}{{2\delta }}\int\limits_0^\pi  {q(\theta )} \cos  \mu \Delta (\theta )\cos  \mu \pi d\theta \\
&  - \frac{{{ \mu ^2}{{a'_2}}\sin b}}{{2\delta }}\int\limits_0^\pi  {q(\theta )} \cos  \mu (2\theta  - \Delta (\theta ))\cos  \mu \pi d\theta  - \frac{{{ \mu ^2}{a_2}\sin b}}{\delta }\sin  \mu \pi \\
&  - \frac{{{ \mu ^2}{{a'_2}}\sin b}}{{2\delta }}\int\limits_0^\pi  {q(\theta )} \sin  \mu \Delta (\theta )\sin  \mu \pi d\theta \\
&  - \frac{{{ \mu ^2}{{a'_2}}\sin b}}{{2\delta }}\int\limits_0^\pi  {q(\theta )} \sin  \mu (2\theta  - \Delta (\theta ))\sin  \mu \pi d\theta \\
&= 0.
\end{array}\]

\noindent Now, let us denote

\begin{equation*}
    \begin{array}{rl}
   P( \mu ,\Delta (\theta ))& =  \frac{1}{2}\int\limits_0^\pi  {q(\theta )} \cos  \mu \Delta (\theta )d\theta\\
     Q( \mu ,\Delta (\theta )) &=  \frac{1}{2}\int\limits_0^\pi  {q(\theta )} \cos  \mu (2\theta  - \Delta (\theta ))d\theta \\
    R( \mu ,\Delta (\theta )) &=    \frac{1}{2}\int\limits_0^\pi  {q(\theta )} \sin  \mu \Delta (\theta )d\theta  \\
    S( \mu ,\Delta (\theta )) &= \frac{1}{2}\int\limits_0^\pi  {q(\theta )} \sin  \mu (2\theta  - \Delta (\theta ))d\theta
    \end{array}.
\end{equation*}

\noindent and we write

\[\begin{array}{rl}
F(\mu ) &=  - \frac{{\mu ^3}{{a'_2}}\sin b}{\delta }\sin \mu \pi  + \frac{{\mu ^2}}{\delta}\left( {{{a'_1}}\sin b - {{a'_2}}\sin bP(\mu ,\Delta (\theta )) }\right.\\
& \,\,\,\,\,\, \left. {- {{a'_2}}\sin bQ(\mu ,\Delta (\theta ))} \right)\cos \mu \pi \\
 &\,\,\,\,\,\, + \frac{\mu ^2}{\delta }\left( { - {a_2}\sin b +  - {a'_2}\sin bR(\mu ,\Delta (\theta )) - {a'_2}\sin bS(\mu ,\Delta (\theta ))} \right)\sin \mu \pi \\
 &\,\,\,\,\,\,+ O\left( {\mu {e^{\vert {Im \mu} \vert \pi }} }\right)\\
\end{array}\]

\section{A Formula For The Regular Trace}

\noindent Let's denote the first term in the last equation of $F(\mu )$ by ${F_0}(\mu )$

\begin{equation}
    {F_0}(\mu ) =  - \frac{{{\mu ^3}{a'_2}}}{\delta }\sin b .\sin \mu \pi. \
\end{equation}

\noindent Let's show the zeros of ${{F_0}(\mu )}$  except the zeros $\mu^{0}_{\mp 0}=\mu^{0}_{\mp 1}=0$ with multiplicity of 4, as follows [similiar to \cite{hira}]

  $$ \mu _n^0=\left\{\begin{array}{rl}
   n-1,&n\geq 1\\
   n+1,&n\leq 1
   \end{array}\right .$$
   
\noindent Here $n \in \mathds{N}\bigcup\{0\} $.
\noindent Let's denote circles of radius $\varepsilon$ with centers at $ \mu _n^0$ points by $\Gamma_{N_{0}}$.From the last equation of $F(\mu )$ and the equation of ${F_0}(\mu )$, on  the contour $\Gamma_{N_{0}}$, we have

\begin{equation}
\begin{array}{l}
\frac{{F(\mu )}}{{{F_0}(\mu )}} = 1 + \frac{1}{\mu }\left\{ \left( { - \frac{{{a'_1}}}{{{a'_2}}} + P(\mu ,\Delta (\theta )) + Q(\mu ,\Delta (\theta ))} \right)\cot \mu \pi  \right.\\
\,\,\,\,\,\,\,\,\,\,\,\,\,\,\,\,\,\,+ \left.  \frac{{{a_2}}}{{{a'_2}}} + R(\mu ,\Delta (\theta )) + S(\mu ,\Delta (\theta )) \right\}\\
\,\,\,\,\,\,\,\,\,\,\,\,\,\,\,\,\,\, + O\left( {\frac{1}{{{\mu. ^2}}}} \right).
\end{array}
\end{equation}

\noindent Expanding $\ln \frac{{F(\mu )}}{{{F_0}(0)}}$ by the Maclaurin formula, we obtain

\[\begin{array}{l}
\ln \frac{{F(\mu )}}{{{F_0}(\mu )}} = \frac{1}{\mu }\left\{ {\left( { - \frac{{{{a'_1}}}}{{{{a'_2}}}} + P(\mu ,\Delta (\theta )) + Q(\mu ,\Delta (\theta ))} \right)\cot \mu \pi }\right.\\
\,\,\,\,\,\,\,\,\,\,\,\,\,\,\,\,\,\,\,\,\,\,\, +\left. { \frac{{{a_2}}}{{{{a'_2}}}} + R(\mu ,\Delta (\theta )) + S(\mu ,\Delta (\theta ))} \right\}\\
\,\,\,\,\,\,\,\,\,\,\,\,\,\,\,\,\,\,\,\,\,\,\, - \frac{1}{{2{\mu ^2}}}{\left( { - \frac{{{{a'_1}}}}{{{{a'_2}}}} + P(\mu ,\Delta (\theta )) + Q(\mu ,\Delta (\theta ))} \right)^2}{\cot ^2}\mu \pi \\
\,\,\,\,\,\,\,\,\,\,\,\,\,\,\,\,\,\,\,\,\,\,\, - \frac{1}{{2{\mu ^2}}}{\left( {\frac{{{a_2}}}{{{{a'_2}}}} + R(\mu ,\Delta (\theta )) + S(\mu ,\Delta (\theta ))} \right)^2}\\
\,\,\,\,\,\,\,\,\,\,\,\,\,\,\,\,\,\,\,\,\,\, + \frac{1}{{{\mu ^2}}}\,\left( { - \frac{{{a'_1}}}{{{{a'_2}}}} + P(\mu ,\Delta (\theta )) + Q(\mu ,\Delta (\theta ))} \right) \\
\,\,\,\,\,\,\,\,\,\,\,\,\,\,\,\,\,\,\,\,\,\, \times\left( {\frac{{{a_2}}}{{{{a'_2}}}} + R(\mu ,\Delta (\theta )) + S(\mu ,\Delta (\theta ))} \right)\cot \mu \pi + O\left( {\frac{1}{{{\mu ^3}}}} \right).
\end{array}\]

\noindent Using the Rouche theorem it follows that ${F(\mu )}$ has the same number of zeros inside the contour as ${{F_0}(0)}$ . Then, we have ${\mu _n} = \mu _n^0 + {\varepsilon _n}$
for sufficiently large n, where $\vert {{\varepsilon _n}} \vert < {\textstyle{\pi  \over 2}}$. Substituting into we get ${\varepsilon _n} = O({\textstyle{1 \over n}})$. We continue making ${\mu _n}$ more precise. Considering the residue theorem, we obtain the asymptotic expression of the eigenvalues as follows:

\[\begin{array}{rl}
{\mu _n} - \mu _n^0 =&   - \frac{1}{{2\pi i}}\oint\limits_{{C_n}} {\ln \frac{{F(\mu )}}{{{F_0}(\mu )}}} d\mu \\
\,\,\,\,\,\,\,\,\,\,\,\,\,\,\,\,\, =&   - \frac{1}{{2\pi i}}\oint\limits_{{C_n}} {\left( { - \frac{{{a'_1}}}{{{a'_2}}} + P(\mu ,\Delta (\theta )) + Q(\mu ,\Delta (\theta ))} \right)} \frac{{\cot \mu \pi }}{\mu }d\mu \\
& \,\,\,\, - \frac{1}{{2\pi i}}\oint\limits_{{C_n}} {\left( {\frac{{{a_2}}}{{{a'_2}}} + R(\mu ,\Delta (\theta )) + S(\mu ,\Delta (\theta ))} \right)} \frac{1}{\mu }d\mu  + O\left( {\frac{1}{{{n^2}}}} \right)\\
\,\,\,\,\,\,\,\,\,\,\,\,\,\,\,\,\, =&   - \frac{1}{{\mu _n^0\pi }}\left\{ { - \frac{{{a'_1}}}{{{a'_2}}} + P(n,\Delta (\theta )) + Q(n,\Delta (\theta ))} \right\} + O\left( {\frac{1}{{{n^2}}}} \right)
\end{array}\]

\noindent Thus, we have proved Theorem 1, which we have stated below.
\begin{theorem} When n approaches infinity, the eigenvalues of the boundary value problem, we have created with (1) differential equation, (2)-(3) boundary conditions and (4)-(5) interface conditions, are expressed by the following asymptotic formula

\begin{equation}
{\mu _n} = \mu _n^0    - \frac{1}{{\mu _n^0\pi }}\left\{ { - \frac{{{a'_1}}}{{{a'_2}}} + P(n,\Delta (\theta )) + Q(n,\Delta (\theta ))} \right\} + O\left( {\frac{1}{{{n^2}}}} \right).
\end{equation}
\end{theorem}
\begin{theorem}

\noindent The regular trace formula for the eigenvalues of the $(1)-(5)$ problem is obtained with the $P(\mu ,\Delta (\theta ))$, $Q(\mu ,\Delta (\theta ))$, $R(\mu ,\Delta (\theta ))$ and $S(\mu ,\Delta (\theta ))$ integrals we have determined, as follows

\[\begin{array}{l}
\mu _{ - 0}^2 + \mu _0^2 + \sum\limits_{n = 1}^\infty  {\left( {\mu _{ - n}^2 + \mu _n^2 - 2{{(n - 1)}^2} + \frac{4}{\pi }\left( { - \frac{{{{a'_1}}}}{{{{a'_2}}}} + P(\mu ,\Delta (\theta )) + Q(\mu ,\Delta (\theta ))} \right)} \right)} \\
\,\,\,\,\,\,\, =  - \frac{2}{\pi }\left( { - \frac{{{{a'_1}}}}{{{{a'_2}}}} + P(0,\Delta (\theta )) + Q(0,\Delta (\theta ))} \right)\\
\,\,\,\,\,\,\,\,\,\,\, - {\left( { - \frac{{{{a'_1}}}}{{{{a'_2}}}} + P(0,\Delta (\theta )) + Q(0,\Delta (\theta ))} \right)^2}\\
\,\,\,\,\,\,\,\,\,\,\, + {\left( {\frac{{{a_2}}}{{{{a'_2}}}} + R(0,\Delta (\theta )) + S(0,\Delta (\theta ))} \right)^2} + O\left( {\frac{1}{{{N_0}}}} \right).\,\,\,\,\,\,\,\,\,\,\,\,\,\,\,\,\,\,\,\,\,\,\,\,\,\,\,\,\,\,\,\,\,\,\,\,\,\,\,\,\,\,\,\,\,\,\,\,\,\,\,\,\,\,\,\,\,\,\,\,\,\,\,\,\,\,\,\,\,\,\,\,\,\,\,(21)
\end{array}\]

\end{theorem}
\begin{proof}

 ${N_0}$ is an integer and denote by ${\Gamma _{{N_0}}}$ the countersclockwise square contours EFGH with $E = ({N_0} - 1 + \varepsilon )(1 - i)$, $F = ({N_0} - 1 + \varepsilon )(1 + i)$, $G = ({N_0} - 1 + \varepsilon )( - 1 + i)$, $H = ({N_0} - 1 + \varepsilon )( - 1 - i)$. Asymptotic formula of ${\mu _n}$ implies that for all sufficiently large ${N_0}$, the numbers  ${\mu _n}$, with $\left| n \right| \le {N_0}$ are inside ${\Gamma _{{N_0}}}$. The numbers  ${\mu _n}$, with $\left| n \right| > {N_0}$ are outside ${\Gamma _{{N_0}}}$. Obviously $\mu _n^0$  do not lie on the contour ${\Gamma _{{N_0}}}$. It follows that

\begin{equation*}
\begin{array}{rl}
\sum\limits_{{\Gamma _{{N_0}}}}^{} {\mu _n^2 + {{(\mu _n^0)}^2}} & =  \mu _{ - 0}^2 + \mu _0^2 + \sum\limits_{n = 1}^{_{{N_0}}} {(\mu _n^2 + \mu _{ - n}^2 - 2{{(n - 1)}^2})} \\
& =  - \frac{1}{{2\pi i}}\oint\limits_{{\Gamma _{{N_0}}}} {2\mu \ln \frac{{F(\mu )}}{{{F_0}(\mu )}}} d\mu \\
& =  - \frac{1}{{2\pi i}}\oint\limits_{{C_n}} {2\left( { - \frac{{{a'_1}}}{{{a'_2}}} + P(\mu ,\Delta (\theta )) + Q(\mu ,\Delta (\theta ))} \right)} \cot \mu \pi d\mu \\
\,\,\,\,\,\,\,\,\,\,\,\,\,\,\,\,\,\,\,\,\,\,\,\,\,\,\,& \,\,\,\, - \frac{1}{{2\pi i}}\oint\limits_{{C_n}} {2\left( {\frac{{{a_2}}}{{{a'_2}}} + R(\mu ,\Delta (\theta )) + S(\mu ,\Delta (\theta ))} \right)} d\mu \\
\,\,\,\,\,\,\,\,\,\,\,\,\,\,\,\,\,\,\,\,\,\,\,\,\,& \,\,\,\,\, + \frac{1}{{2\pi i}}{\oint\limits_{{C_n}} {\left( { - \frac{{{a'_1}}}{{{a'_2}}} + P(\mu ,\Delta (\theta )) + Q(\mu ,\Delta (\theta ))} \right)} ^2}\frac{{{{\cot }^2}\mu \pi }}{\mu }d\mu \\
\,\,\,\,\,\,\,\,\,\,\,\,\,\,\,\,\,\,\,\,\,\,\,\,\,& \,\,\,\,\, + \frac{1}{{2\pi i}}{\oint\limits_{{C_n}} {\left( {\frac{{{a_2}}}{{{a'_2}}} + R(\mu ,\Delta (\theta )) + S(\mu ,\Delta (\theta ))} \right)} ^2}\frac{1}{\mu }d\mu \\
\,\,\,\,\,\,\,\,\,\,\,\,\,\,\,\,\,\,\,\,\,\,\,\,\,\,& \,\,\,\, + \frac{1}{{2\pi i}}\oint\limits_{{C_n}} {2\left( { - \frac{{{a'_1}}}{{{a'_2}}} + P(\mu ,\Delta (\theta )) + Q(\mu ,\Delta (\theta ))} \right)} \\
\,\,\,\,\,\,\,\,\,\,\,\,\,\,\,\,\,\,\,\,\,\,\,\,\,& \,\,\,\,\, \times \left( {\frac{{{a_2}}}{{{a'_2}}} + R(\mu ,\Delta (\theta )) + S(\mu ,\Delta (\theta ))} \right)\frac{{\cot \mu \pi }}{\mu }d\mu  + O\left( {\frac{1}{{{N_0}}}} \right).
\end{array}    
\end{equation*}

\noindent And we get

\[\begin{array}{l}
\mu _{ - 0}^2 + \mu _0^2 + \sum\limits_{n = 1}^{_{{N_0}}} {(\mu _n^2 + \mu _{ - n}^2 - 2{{(n - 1)}^2})} \\
\,\,\,\,\,\,\,\,\,\,\,\,\,\,\,\,\,\,\,\,\,\,\,\,\,\,\,\,\, =  - 2\left( { - \frac{{{a'_1}}}{{{a'_2}}} + P(n,\Delta (\theta )) + Q(n,\Delta (\theta ))} \right)\frac{{(2({N_0} - 1) + 1)}}{\pi }\\
\,\,\,\,\,\,\,\,\,\,\,\,\,\,\,\,\,\,\,\,\,\,\,\,\,\,\,\,\,\,\, - \frac{2}{\pi }\left( { - \frac{{{a'_1}}}{{{a'_2}}} + P(0,\Delta (\theta )) + Q(0,\Delta (\theta ))} \right)\\
\,\,\,\,\,\,\,\,\,\,\,\,\,\,\,\,\,\,\,\,\,\,\,\,\,\,\,\,\,\, - {\left( { - \frac{{{a'_1}}}{{{a'_2}}} + P(0,\Delta (\theta )) + Q(0,\Delta (\theta ))} \right)^2}\\
\,\,\,\,\,\,\,\,\,\,\,\,\,\,\,\,\,\,\,\,\,\,\,\,\,\,\,\,\,\, + {\left( {\frac{{{a_2}}}{{{a'_2}}} + R(0,\Delta (\theta )) + S(0,\Delta (\theta ))} \right)^2} + O\left( {\frac{1}{{{N_0}}}} \right).
\end{array}\]
\noindent again by residue calculation. So, this last equation implies the following equation

\[\begin{array}{l}
\mu _{ - 0}^2 + \mu _0^2 + \sum\limits _{n = 1}^\infty  {\left( {\mu _{ - n}^2 + \mu _n^2 - 2{{(n - 1)}^2} + \frac{4}{\pi }\left( { - \frac{{{a'_1}}}{{{a'_2}}} + P(\mu ,\Delta (\theta )) + Q(\mu ,\Delta (\theta ))} \right)} \right)} \\
 =  - \frac{2}{\pi }\left( { - \frac{{{a'_1}}}{{{a'_2}}} + P(0,\Delta (\theta )) + Q(0,\Delta (\theta ))} \right)\\
\,\,\,\, - {\left( { - \frac{{{a'_1}}}{{{a'_2}}} + P(0,\Delta (\theta )) + Q(0,\Delta (\theta ))} \right)^2}\\
\,\,\,\, + {\left( {\frac{{{a_2}}}{{{a'_2}}} + R(0,\Delta (\theta )) + S(0,\Delta (\theta ))} \right)^2} + O\left( {\frac{1}{{{N_0}}}} \right).
\end{array}\]

\noindent Thus, by approximating ${N_0}$ to  $\infty $, we get the regular trace formula (21) which we want, and our proof is complete.
\end{proof}

\end{document}